\documentclass[12pt,reqno]{amsart}
\usepackage{amsmath}
\usepackage{amssymb}
\usepackage{amsfonts}
\numberwithin{equation}{section}
\newcommand{\D}{\mathbb{D}\,}
\newcommand{\T}{\mathbb{T}\,}
\newcommand{\Z}{\mathbb{Z}\,}
\newcommand{\C}{\mathbb{C}\,}
\newcommand{\NN}{\mathbb{N}\,}

\newcommand{\HH}{{\mathcal H}}
\newcommand{\E}{{\mathcal E}}
\newcommand{\M}{{\mathcal M}}
\newcommand{\U}{{\mathcal U}}
\newcommand{\A}{{\mathcal A}}
\newcommand{\e}{{\varepsilon}}
\newcommand{\SU}{{\mathcal S\mathcal U}}

\DeclareMathOperator{\dist}{dist}
\DeclareMathOperator{\card}{Card}
\DeclareMathOperator{\supp}{supp}

\newtheorem{th1}{{\bf Theorem}}[section]

\newtheorem{thm}[th1]{{\bf Theorem}}
\newtheorem{lem}[th1]{{\bf Lemma}}

\theoremstyle{definition}
\newtheorem{rem}[th1]{{\bf Remark}}
\newtheorem{que}[th1]{{\bf Question}}
\newtheorem{exemp}[th1]{{\bf Example}}

\newcommand{\lt}{\log_2}

\newcommand{\ts}{\tilde{\Sigma}}

\begin{document}
\author[A. Borichev, Yu. Lyubarskii]
{Alexander Borichev, Yurii Lyubarskii}

\title[Uniqueness theorems for Korenblum type spaces]
{Uniqueness theorems for Korenblum type spaces}

\begin{abstract}
For a scale of spaces $X$ of functions analytic in the unit disc,
including the Korenblum space, and for some natural families $\mathcal E$
of uniqueness subsets for $X$, we describe minorants for $(X,\mathcal E)$, that is
non-decreasing functions $M:(0,1)\to(0,\infty)$ such that $f\in X$, $E\in\mathcal E$,
and $\log|f(z)|\le -M(|z|)$ on $E$ imply $f=0$. We give an application of this result
to approximation by simple fractions with restrictions on the coefficients.
\end{abstract}

\date{\today}

\maketitle

\section{Introduction.}

Given a topological space $X$ of analytic functions in the unit disc $\D$
and a class $\E$ of subsets $E$ of $\D$, we call a non-decreasing positive continuous function
$M:[0,1)\to(0,\infty)$ a minorant for the pair $(X,\E)$ and write $M\in\M(X,\E)$ if
for every $f\in X$ and $E\in\E$ such that
\begin{equation}
\log|f(z)|\le -M(|z|),\qquad z\in E,
\label{one}
\end{equation}
we have $f=0$.

Clearly, $\M(X,\E)\ne\emptyset$ implies that $\E\subset \U(X)$,
where $\U(X)$ is the class of all uniqueness subsets $E$ for the
space $X$: $E\in\U(X)$ if and only if $E\subset\D$, $E$ has no
limit points in $\D$, and
$$
f\in X,\quad f{\bigm|}E=0\quad
$$
imply that $f=0$.

Suppose that $H^\infty\subset X$ and $X$ is not too large, that is $X\subset A(\lambda)$, for some function $\lambda:[0,1)\to[0,+\infty)$, where
$$
A(\lambda)=\bigl\{f\in\text{Hol}(\D):\sup_{z\in\D}|f(z)|/\lambda(|z|)<\infty\bigr\}.
$$
Then the family $\M(X,\U(X))$ is empty. The reason is that the class $\U(X)$
contains subsets $E\subset\D$ which are "too concentrated": $E$ may be the union
of clusters $E_j$ of nearby points in such a way that the estimate \eqref{one}
at one point $x\in E_j$ implies a similar estimate (eventually, with $M$ replaced by $cM$)
on the whole $E_j$. That is why we need to consider only elements in the family
$\U(X)$ which are ``sufficiently separated''.

In \cite{borichev3}, Lyubarskii and Seip
deal with the case $X=H^\infty$. They consider the class $\SU(H^\infty)$ of hyperbolically
separated subsets $E$ of $\D$ that are uniqueness subsets for $H^\infty$ and prove that
$$
M\in \M(H^\infty,\SU(H^\infty)) \iff \int_0\frac{dt}{tM(1-t)}<\infty.
$$
For other uniqueness theorems of such type for the space
$H^\infty$ see \cite{Us} \cite{Da}, \cite{Ha}, \cite{borichev4}
\cite{Va}.
 We also
refer the reader to the survey article \cite{ee} and the
references therein.

In the present article we work with the scale of spaces
\begin{gather*}
\A^r_s=\bigl\{f\in\text{Hol}(\D): \log|f(z)|\le r\log^s\frac 1{1-|z|}+c_f\bigr\},\qquad r,s>0,\\
\A_s=\bigcup_{r<\infty}\A^r_s.
\end{gather*}
We have $H^\infty\subset\A_s\subset\A_1\subset\A_t$, $0<s<1<t$, where $\A_1$ is the
{\em Korenblum space},
$$
\A_1=\bigl\{f\in\text{Hol}(\D): |f(z)|\le \frac {c_f}{(1-|z|)^{c'_f}}\bigr\}.
$$

The uniqueness subsets for $\A_s$ are described by Korenblum \cite{borichev2}
(1975, $s=1$), Pascuas \cite{Pas} (1988, $0<s<1$), and Seip \cite{borichev5} (1995, $s>0$). For $0<s<1$, we define
$\SU(\A_s)$ as the class of hyperbolically separated subsets $E$ of $\D$
that are uniqueness subsets for $\A_s$.

\begin{thm} For $0<s<1$,
$$
M\in \M(\A_s,\SU(\A_s)) \iff \int_0\frac{dt}{tM(1-t)}<\infty.
$$
\label{te1}
\end{thm}

The case $s=1$ should be treated differently because $\U(\A_1)$
contains no hyperbolically separated subsets of $\D$.
However, for each $r>0$, the class $\U(\A^r_1)$ contains hyperbolically
separated subsets of $\D$
(see \cite{Ho}, \cite{borichev2}, \cite{borichev5}, and \cite{hkkbook}).
That is why we define $\SU(\A_1)$ as the class of $E\subset\D$ satisfying the following property: for
every $r<\infty$ there exists a
hyperbolically separated subset $E_r$ of $E$ such that $E_r\in\U(\A^r_1)$.

\begin{thm}
$$
M\in \M(\A_1,\SU(\A_1)) \iff \int_0\frac{dt}{tM(1-t)}<\infty.
$$
\label{te2}
\end{thm}

For $s>0$, we set $\beta=\max(0,s-1)$, and introduce
$$
\rho_s(z)=(1-|z|)\Bigl(\log\frac1{1-|z|}\Bigr)^{-\beta/2};
$$
this is a natural measure of separateness for points from uniqueness sets for $\A_s$,
see Section~\ref{sect3} below for more details. Given $z\in\C$, $a>0$, we denote
$D_z(a)=\{\zeta\in\C:|z-\zeta|<a\}$. When $s$ is fixed, we denote $D^t_z=D_z(t\rho_s(z))$.
We say that $E\subset\D$ is $s$-{\em separated} (with constant $\e>0$) if
$$
D^\e_\lambda \cap D^\e_\mu=\emptyset,\qquad \lambda,\mu\in E,\quad \lambda\ne\mu.
$$
For $0<s\le 1$, this notion coincides with the usual hyperbolic separateness.

We define $\SU(\A_s)$, $s>1$, as the class of $E\subset\D$
such that for every $r<\infty$ there exists an
$s$-separated subset $E_r$ of $E$ such that $E_r\in \U(\A^r_s)$.
\medskip

\begin{thm} For $s>1$,
$$
M\in \M(\A_s,\SU(\A_s)) \iff \int_0\Bigl(\log \frac 1t\Bigr)^{s-1}\frac{dt}{tM(1-t)} <\infty.
$$
\label{te3}
\end{thm}

\begin{rem} The uniqueness results we formulate here involve three
factors: possible growth of the functions in the function spaces $\A^s$,
the separateness of the uniqueness sets $E$ (determined by the distance function $\rho_s$),
and the integral condition on the minorant.
Surprisingly, it turns out that the size of the minorant for which uniqueness holds
depends more upon the measure of separateness of subsets $E$ than
on the growth of functions in $\A_s$. In fact, in Theorems~\ref{te1}--\ref{te3},
when the integrals diverge, we can find $E\in\SU(\A_s)$ and a \textit{bounded}
$f\ne0$ holomorphic in $\D$ satisfying the estimate \eqref{one}, see Section~\ref{sect2}.
Furthermore, the classes $\SU(\A_s)$ are different for different values of $s$.
In Section~\ref{sect7} we give, for each $s>1$, an example of a set in $\SU(\A_s)$
which contains no hyperbolically separated non-Blaschke sequence. 
\end{rem}

\begin{rem}  
The proofs of our results demand new ideas and techniques in addition to
those already used in the $H^\infty$ case. In particular, in the standard dyadic
decomposition of the unit disc, we need to deal with large groups of adjacent squares, 
see for example Fig. 1 below.
\end{rem}

\begin{que} How to get analogous results for the Bergman space (no description
of uniqueness subsets is known yet), for the spaces $\A^r_s$, $0<s<1$?
\end{que}

The article is organized as follows. In Section~\ref{sect2} we formulate
two theorems that imply Theorems~\ref{te1}--\ref{te3}.
To prove these results we use some information on uniqueness subsets
for the classes $\A^r_s$, $\A_s$ obtained by Korenblum and Seip.
This information is contained in Section~\ref{sect3}.
Section~\ref{sect4} contains the proof of our sufficient condition
on the minorants, and Section~\ref{sect5} is devoted to the construction of
an example which shows that this condition is also necessary.
As an application of our results, in section~\ref{sect6} we prove a
statement on approximation by simple fractions with restrictions on the coefficients
in the space of  functions analytic in the unit disc and infinitely smooth
up to the boundary. Final remarks are given in Section~\ref{sect7}. We
compare our results to those by Pau and Thomas \cite{borichev4}, and give an example
illustrating an important difference between the classes $\SU(\A_s)$, $s\le 1$, and $\SU(\A_s)$, $s>1$.

The authors are thankful to Nikolai Nikolski and Pascal Thomas for 
valuable comments. In particular, the definition of the class $\SU(\A_1)$
was first formulated during our discussions with Pascal Thomas.
\medskip

\section{Main results}
\label{sect2}

In this section we first introduce (in a way independent of $s>0$) a family $\mathfrak A_s$ of
``massive'' subsets $E$ of the unit disc on which our uniqueness theorem holds for large $M$.
After that, we formulate our main theorems. They imply Theorems~\ref{te1}--\ref{te3},
and show that for small $M$ we can find bounded $f$ satisfying $\log|f|\le-M$ on
uniqueness subsets $E$.

Given $s>0$, we define a class $\mathfrak A_s$ of subsets of the unit disc
in the following way: $E\in\mathfrak A_s$ if and only if for every $r<\infty$
there exists an $s$-separated subset $E_r$ of $E$ such that $E_r\in \U(\A^r_s)$. Then
$$
\SU(\A_s)\subsetneqq \mathfrak A_s,\quad s<1,\qquad \SU(\A_s)=\mathfrak A_s,\quad s\ge 1.
$$

Furthermore, we assume that $M:(0,1)\to(0,\infty)$ is a non-decreasing continuous function.
Recall that $\beta=\max(0,s-1)$.

\begin{thm} Let $s>0$, $E\in\mathfrak A_s$, and let
\begin{equation}
\int_0 \frac 1t\Bigl(\log\frac 1t \Bigr)^\beta\frac{dt}{M(1-t)}<\infty.
\label{convergence}
\end{equation}
If $f\in \A_s$ and
\begin{equation}
\log|f(\lambda)|\le -M(|\lambda|), \qquad  \lambda\in E,
\label{estiamte}
\end{equation}
then $f=0$.
\label{thm1}
\end{thm}

\begin{thm} Let $s>0$, and let
\begin{equation}
\int_0 \frac 1t \Bigl(\log\frac 1 t \Bigr)^\beta \frac {dt}{M(1-t)}=\infty.
\label{divergence}
\end{equation}
Then there exist $E\in \SU(\A_s)$ and $f\in H^\infty\setminus\{0\}$ such that
$$%\begin{equation}
\log |f(\lambda)|\le-M(|\lambda|),\qquad \lambda\in E.
%\label{mainestiamte}
$$%\end{equation}
\label{thm2}
\end{thm}

We end this section by mentioning the following analog of Theorem~\ref{thm1} for the space $H^\infty$:

\begin{rem} If $s>0$, $E\subset\D$ is $s$-separated and does not satisfy the Blaschke condition, and if
\eqref{convergence} holds, then any $f\in H^\infty$ satisfying \eqref{estiamte} is equal to $0$, see
\cite[Corollaries~4.7, 4.8]{ee}.
\end{rem}
\medskip

\section{Uniqueness subsets for $\A_s$.}
\label{sect3}

In this section we formulate a (known) uniqueness theorem for spaces $\A^r_s$.

First we need notion of $s$-{\em entropy}, $s>0$. Given a finite set $F\subset \mathbb{T}$,
let $\{I_j\}$, $I_j=\{e^{i\theta}: |\theta-\theta_j|<t_j\}$ be
the set of all disjoint arcs complementary to $F$. We define the $s$-entropy of $F$ as
$$
\varkappa_s(F)=\sum|I_j|\Bigl(\log\frac e{|I_j|}\Bigr)^s,
$$
where $|I_j|$ is the normalized length of the arc $I_j$, i.e. $t_j/\pi$.

Consider the {\em Korenblum stars} associated with the set $F$:
$$
\mathcal{K}^t_1(F)=\bigl\{z\in \D:\dist(z,\mathbb{T})>t\dist(z,F)\bigr\},\qquad 0<t<1,
$$
where $\dist$ stays for usual Euclidian distance.

For every sequence $E\subset \D$ and a finite set $F\subset \mathbb{T}$, we define
$$
\Sigma_t(F,E)=\sum_{\lambda\in E\cap\mathcal{K}^t_1(F)}(1-|\lambda|)
$$
and
$$
\delta_s(E)=\limsup_{\varkappa_s(F)\to\infty}\frac{\Sigma_t(F,E)}{\varkappa_s(F)}.
$$
The value $\delta_s(E)$ does not depend on $t\in(0,1)$, see, for example, the dicussion in
\cite[Section~4.2]{hkkbook}.

\begin{thm} {\rm(\cite{borichev2, Pas, borichev5})} For a sequence $E$ to be a set of uniqueness
for $\A^r_s$ it is necessary that $\delta_s(E)\ge r$ and sufficient
that $\delta_s(E)>r$.
\label{korteo}
\end{thm}

For $s=1$ this result is proved in \cite{borichev2}, for $0<s<1$ it is proved in
\cite{Pas}, for $s>0$ it is formulated in \cite{BM, borichev5}. To get the proof in the case
$s>0$, one can just follow the argument in \cite{borichev2} or in \cite[Sections 4.4, 4.5]{hkkbook}.
The only modification to be done is to replace the entropy $\sum_j\mathbf k_1(|I_j|)$, $\mathbf k_1(x)=x\log^+(1/x)$,
for the case $s=1$ by the entropy $\sum_j \mathbf k_s(|I_j|)$, $\mathbf k_s(x)=x[\log^+(C_s/x)]^s$,
for the case $s>1$, with large $C_s$ and $\mathbf k_s$ concave on $(0,1)$. 

In the case $s\in(0,1)$ the necessary conditions and the sufficient conditions merge, see \cite[Theorem~1.1]{borichev5}).

Denote by $\mathcal{Q}$ a dyadic family of Whitney squares:
$$
Q_{n,k}=\bigl\{ re^{i\theta}\in \D: \theta\in (\pi k2^{-n},\pi(k+1)2^{-n}),
\, 2^{-n-1}<1-r<2^{-n}\bigr\},
$$
with $n\ge 0$, $0\le k<2^{n+1}$. We denote by
$|Q|$ the diameter of $Q\in \mathcal{Q}$.
Furthermore, we define $\mathcal{K}_1(F)=\mathcal{K}^{1/10}_1(F)$, and
\begin{gather*}
\tilde{\Sigma}(F,E)=\sum_{Q\in \mathcal Q:\,Q\subset \mathcal{K}_1(F)}\card(Q\cap E)|Q|,\\
\tilde{\delta}_s(E)=\limsup _{\varkappa_s(F)\to\infty}\frac{{\tilde{\Sigma}}(F,E)}{\varkappa_s(F)}.
\end{gather*}

For each finite $F\subset \mathbb{T}$ we have
$$
\frac 12 \Sigma_{9/10}(F,E)\le\tilde{\Sigma}(F,E)\le 2 \Sigma_{1/10}(F,E).
$$
Hence, for a sequence $E$ to be a set of uniqueness for $\A^r_s$ it is necessary that
${\tilde{\delta}}_s(E)\ge r/2$ and sufficient that  ${\tilde{\delta}}_s(E)>2r$.
The discrepancy in the right hand sides in these conditions does not play any role
since we deal with the whole algebra $\A_s$.

Now we obtain that, for each $s\ge 1$ and $r$, there exists
$E\subset\D$ which is a set of uniqueness for $A^r_s$ and such that
$$
\card(Q\cap E)\le C\Bigl(\frac 1{\log |Q|}\Bigr)^{s-1},\qquad Q\in \mathcal Q,
$$
in particular $E$ may be an $s$-separated set.

These results motivate our definitions of classes $\SU(\A_s)$:
when $s\le 1$ the spaces $A_s^r$ admit sets of uniqueness which are hyperbolically separated.
On the other hand, when $s\ge 1$, every hyperbolically separated
(or even just $s$-separated) set is a set of non-uniqueness for $A_s$.
\medskip

\section{Proof of Theorem~\ref{thm1}}
\label{sect4}

{\bf (A)} Here we establish a helpful asymptotical relation for $M$.
The property \eqref{convergence} implies that
\begin{equation}
\Bigl(\log \frac 1t\Bigr)^{\beta+1}=o(M(1-t)),\qquad t\to 0.
\label{regt}
\end{equation}
Otherwise, we could find $\e>0$, $x_n>0$, $x_n\to 0$, such that
$$
\Bigl(\log\frac 1{x_n}\Bigr)^{\beta+1}\ge \e M(1-x_n).
$$
Since $M$ is monotonic, for small $y>0$ and for every $x_n<y$ we obtain
\begin{gather*}
\frac \e{2(\beta+1)}\ge\int^y_0\frac 1t\Bigl(\log \frac 1t\Bigr)^\beta \frac{dt}{M(1-t)}\\
\ge \int^y_{x_n}\frac 1t\Bigl(\log \frac 1t\Bigr)^\beta
\cdot\e\Bigl(\log \frac 1{x_n}\Bigr)^{-\beta-1}dt\\
=\frac \e{\beta+1}\Bigl(\log \frac 1{x_n}\Bigr)^{-\beta-1}
\Bigl[\Bigl(\log \frac 1{x_n}\Bigr)^{\beta+1}-\Bigl(\log \frac 1{y}\Bigr)^{\beta+1}\Bigr]
\end{gather*}
which is impossible for sufficiently small $x_n$.
\medskip

{\bf (B)} From now on we fix $E$, $M$, and $f\in \A^r_s$ satisfying the hypothesis of Theorem~\ref{thm1}.
Here we obtain estimates on the set where $f$ is small.
We may assume that $f(0)\ne 0$,
\begin{equation}
\log|f(z)|\le r\left(\log\frac 1{1-|z|}\right)^s, \qquad z\in\D.
\label{festimate}
\end{equation}
Denote by $Z_f$ the sequence of the zeros of $f$ taking into account
their multiplicities. Since ${\tilde{\delta}}_s(Z_f)\le 2r$, we can find
small $\eta>0$ and $E'\subset E$ such that
\begin{gather}
{\tilde{\delta}}_s(E')>0, \label{Mdensity}\\
D_\mu^{2\eta} \cap D_\lambda^{2\eta}=\emptyset,\qquad \lambda,\mu\in E',
\quad \lambda \ne \mu,         \label{2etaseparation}\\
D_\lambda^{2\eta}\cap Z_f =\emptyset, \qquad \lambda\in E'. \label{Mzero}
\end{gather}

Furthermore, in the case $s>1$ we need to consider
\begin{equation}
\mathcal{Q}_f=\Bigl\{Q\in \mathcal{Q}: \card(Q\cap Z_f)>
\Bigl(\log\frac 1{d_Q}\Bigr)^{s-1/2}\Bigr\}.                         \label{badQ}
\end{equation}
We denote
$$
E'_f=E'\cap(\cup_{Q\in \mathcal{Q}_f}).
$$
It follows from \eqref{2etaseparation} and \eqref{badQ} that for every $Q\in \mathcal{Q}_f$,
$$
\card(Q\cap E')\le C \card(Q\cap Z_f) \Bigl(\log\frac 1{d_Q}\Bigr)^{-1/2}
$$
and, since ${\tilde{\delta}}_s(Z_f)\le 2r$, we obtain ${\tilde{\delta}}_s(E'_f)=0$.
Thus, we can replace $E'$ by $E'\setminus E'_f$ in relations
\eqref{Mdensity}, \eqref{2etaseparation}, and \eqref{Mzero}. In what
follows we assume that
\begin{equation}
E'\cap(\cup_{Q\in \mathcal{Q}_f}Q)=\emptyset.                     \label{goodlambdas}
\end{equation}

Fix $\mu\in E'$. Since $\eta$ is small, the function $u_\mu$,
$$
u_\mu(z)=\log |f(z)|-2r\Bigl(\log\frac 1{1-|\mu|}\Bigr)^s
$$
is harmonic and non-positive on $D_\mu^{2\eta}$. Furthermore,
$$
u_\mu(\mu)\le -M(|\mu|).
$$
Applying Harnack's inequality to $u_\mu$, for some $C>0$ we obtain
$$
u_\mu(z)\le -C M(|\mu|), \qquad z\in  D_\mu^\eta.
$$
Hence, by \eqref{regt}, for small $\eta>0$,
$$
\log|f(z)|\le -C M(|\mu|), \qquad z\in D_\mu^\eta,
$$
with $C>0$ independent of $\mu\in E'$.

Denoting  $E_\eta =\cup_{\mu\in E'}D_\mu^\eta$ and using that
$1-|z|<2(1-|\mu|)$, $z\in D_\mu^\eta$, we conclude that
$$
\log|f(z)|\le -CT(2(1-|z|)),\qquad z\in E_\eta,
$$
where we use the notation $T(r)=M(1-r)$, $0<r\le 1$.

Next, we need the notion of $\alpha$-area. For
$\Omega\subset \mathbb{D}$, $\alpha\ge 0$, we set
$$
\HH_\alpha(\Omega)=\int_\Omega\frac 1{1-|z|}\Bigl(\log\frac 1{1-|z|}\Bigr)^{\alpha} dm_2(z),
$$
where $dm_2$ is area Lebesgue measure.
Then the limit
$$
\limsup_{\varkappa_s(F)\to\infty}\frac{\HH_\beta(E_\eta \cap \mathcal K^t_1(F))}
{\varkappa_s(F)}
$$
does not depend on $t\in(0,1)$, see again the dicussion in
\cite[Section~4.2]{hkkbook}.
Hence, by \eqref{Mdensity},
\begin{equation}
\limsup_{\varkappa_s(F)\to\infty}\frac{\HH_\beta(E_\eta \cap \mathcal{K}_1(F))}
{\varkappa_s(F)}=
\limsup_{\varkappa_s(F)\to\infty}\frac{\ts(F,E')}{\varkappa_s(F)}>0.
\label{bigarea}
\end{equation}
\medskip

{\bf (C)} Next we introduce an auxiliary subdomain $\mathcal{K}_2(F)$ of the unit disc.
Take any smo\-oth function  $\phi$  on $[-1,1]$ satisfying the relation
$$
\frac {1-|x|^2}{200} \le \phi(x)\le   \frac {1-|x|^2}{100}, \qquad x\in[-1,1].
$$
For each finite $F\in \mathbb{T}$ we define
$$
\mathcal{K}_2(F)=\bigcup_j \Bigl\{re^{i\theta}:|\theta-\theta_j|\leq t_j,
0\le r\le 1-t_j^2\phi \Bigl(\frac{\theta-\theta_j}{t_j}\Bigr)\Bigr\}.
$$

We consider the conformal mapping $w_F: \mathbb{D} \to \mathcal{K}_2(F)$, determined
by the conditions $w_F(0)=0$ and $w_F'(0)>0$. Then, for some positive $c$ independent of $F$,
we have
\begin{gather}
\frac 1 c \leq |w'_F(z)| \le c, \qquad z\in \D,        \label{derivative}  \\
\frac 1 c \le \frac {1-|w_f(z)|}{1-|z|} \le c, \qquad z\in w^{-1}_F(\mathcal{K}_1(f)),
\label{distortion}          \\
\int_\T\Bigl(\log\frac 1{1-|w_f(z)|}\Bigr)^s dm(z) \le c \varkappa_s(F),  \label{nevestimate}
\end{gather}
where $dm$ is normalized Lebesgue measure.
To verify these assertions we use the Kellogg theorem on the boundary regularity
of the conformal mapping, see \cite[Chapter 10, Section 1, Theorem 6]{gol}.
Since $\partial \mathcal{K}_2(F)$ is Lipschitz smooth uniformly in $F$, we obtain
that $\log w_f'$ is continuous in the closed unit disc and, for some $C$ independent of $F$,
$$
|\log w'_F(z)| \le C, \qquad z\in \mathbb{D}.
$$
Hence
\begin{gather*}
e^{-C}\le |w_F'(z)|\le e^{C}, \qquad z\in\D,\\
e^{-C}\le |(w_F^{-1})'(z)|\le e^{C}, \qquad z\in \mathcal{K}_2(F).
\end{gather*}
This gives (\ref{nevestimate}) and (\ref{derivative}). Furthermore, for
$z\in \mathcal{K}_1(F)$,
$$
1-|z| \asymp \dist(z, \partial \mathcal{K}_2(F)),
$$
and the Koebe distortion theorem (see, for example, \cite[Section 1.3, Corollary 1.4]{pomm})
gives (\ref{distortion}).
\medskip

{\bf (D)} In this part of the proof we verify that $f$ belongs to the Nevanlinna class in the domain $\mathcal{K}_2(F)$
and use this fact to evaluate the size of the set of small values of $f$.
We define the function
$$
g=f\circ w_F.
$$
The function $g$ extends continuously to $\T\setminus w_F^{-1}(F)$.
Furthermore, the estimates \eqref{festimate} and \eqref{derivative} show that
if $\zeta\in w_F^{-1}(F)$, then for some $N$, the function $z\mapsto (z-\zeta)^Ng(z)$
is bounded in a neighborhood of $\zeta$. Hence, $g$ belongs to the Nevanlinna class.
Combining (\ref{festimate}) and (\ref{nevestimate}) we obtain
$$
\int_\mathbb{T}\log^+|g(\zeta)|dm(\zeta) \le C \varkappa_s(F),
$$
with $C$ independent of $F$.
Denote
$E^*_\eta =w_F^{-1}(E_\eta\cap \mathcal{K}_1(F))$.
Then
$$
\log|g(z)|\le -C T(A(1-|z|)),\qquad z\in E^*_\eta,
$$
with some constant $A$ independent of $F$.

Combining \eqref{bigarea} and \eqref{distortion} we obtain
\begin{equation}
\limsup_{\varkappa_s(F)\to\infty}\frac{\HH_\beta(E^*_\eta)}{\varkappa_s(F)}\ge c >0.
\label{newbigarea}
\end{equation}
We use the Nevanlinna factorization theorem (see, e.g. \cite{factorization})
to represent $\log|g|$ in the form
$$
\log|g(z)|=(\mathcal{P}*\mu)(z)+\log |B(z)|,
$$
where $\mathcal{P}$ is the standard Poisson kernel,
$\mu$ is a finite measure on $\mathbb{T}$ and
$$
B(z)=\prod_{\zeta \in Z_g} b_\zeta(z), \qquad
b_\zeta(z)=\frac{\zeta-z}{1-\bar{\zeta}z}\frac {\bar{\zeta}}{|\zeta|},
$$
is the Blaschke product corresponding to the zero sequence $Z_g$ of $g$.

Then $\mu=\mu^+-\mu^-$, where $\mu^+$, $\mu^-$ are non-negative measures and
$$
\mu^-(\T)+\log|g(0)|\le \mu^+(\T)\le \int_\T\log^+|g(\zeta)|dm(\zeta) \le C \varkappa_s(F).
$$
Furthermore,
$$
\log|g(0)|+\sum_{\zeta\in Z_g}\log \frac 1{|\zeta|} \le C \varkappa_s(F).
$$
It follows from \eqref{Mzero} and \eqref{distortion} that, for some small positive $\e$
independent of $F$,
\begin{equation}
D_\zeta^\e\cap E^*_\eta=\emptyset,\qquad \zeta\in Z_g.   \label{d15}
\end{equation}
To each zero $\zeta=re^{i\theta}\in Z_g$ we associate the arc
$I_\zeta=\{e^{it}: |t-\theta|\leq 1-r\}$.

It is an easy (and, probably, well-known fact) that
\begin{equation}
\log|b_\zeta(z)|\ge -c(\mathcal{P}*\chi_{I_\zeta})(z)
\label{blaschkefactorestimate}
\end{equation}
for $|z-\zeta|\ge\e(1-|z|)$, with $c$ depending only on $\e$.
Indeed, since both functions are harmonic in the domain
$\Omega=\D\setminus \overline{D_z(\e(1-|z|))}$ and continuous on $\bar\Omega\setminus\partial I_\zeta$,
it is sufficient to compare them on $\partial\Omega$,
which is elementary.

Thus, for $s\le 1$ we obtain
\begin{gather*}
\log|B(z)|\ge -C(\mathcal{P}*\nu_0)(z), \qquad z\in E^*_\eta,\\
\nu_0(\T)\le C \varkappa_s(F),
\end{gather*}
where
\begin{equation}
d\nu_0=\sum_{\zeta\in Z_g}\chi_{I_\zeta}dm.
\label{nunol}
\end{equation}
\medskip

{\bf (E)} For $s>1$ we need an additional estimate for points of $Z_g$ close to $E^*_\eta$.
Given a point $z_0\in E^*_\eta$ and small $\e>0$, we define
$$
Z^{(1)}_g=\{\zeta \in Z_g: |\zeta-z_0|<\e(1-|z_0|)\}, \qquad Z^{(2)}_g=Z_g\setminus Z^{(1)}_g,
$$
and
$$
B_1(z)=\prod_{\zeta \in Z_g^{(1)}} b_\zeta(z), \qquad
B_2(z)=\prod_{\zeta \in Z_g^{(2)}} b_\zeta(z).
$$
Then $B=B_1B_2$.

The estimate \eqref{blaschkefactorestimate} gives
$$
\log|B_2(z_0)|\geq - C (\mathcal{P}*\nu_0)(z_0),
$$
with $\nu_0$ defined by \eqref{nunol}.
Next we choose $z'_0$ such that $|z_0-z'_0|=2\e(1-|z'_0|)$.
Combining \eqref{blaschkefactorestimate} and Harnack's inequality we obtain
$$
\log|B_1(z'_0)|\ge-C(\mathcal{P}*\nu_0)(z'_0)
\ge-C_1(\mathcal{P}*\nu_0)(z_0).
$$
It remains to compare $|B_1(z'_0)|$ and $|B_1(z_0)|$.

By \eqref{d15}, for each $\zeta \in Z^{(1)}_g$ we have
$$
\e(1-|\zeta|)\Bigl(\log\frac 1{1-|\zeta|}\Bigr)^{(1-s)/2}\le |\zeta-z_0|<\e(1-|z_0|),
$$
and
$$
\e(1-|z_0|)<|\zeta-z'_0|<3\e(1-|z_0|).
$$
A direct estimate gives now
$$
\bigl|\log|b_\zeta(z_0)|-\log|b_\zeta(z'_0)|\bigr|\le C\log \log \frac 1 {1-|z_0|}.
$$
On the other hand, it follows from (\ref{goodlambdas}) that
$$
\card Z_g^{(1)} \le C \Bigl(\log\frac 1 {1-|z_0|}\Bigr)^{s-1/2}.
$$
Therefore
\begin{multline*}
\log|B_1(z_0)|\\ \ge -C(\mathcal{P}*\nu_0)(z_0)
- C_1 \Bigl( \log \frac 1 {1-|z_0|} \Bigr)^{s-1/2}\log\log \frac 1 {1-|z_0|},
\quad z_0\in E^*_\eta.
\end{multline*}

Finally, for any $s>0$ we set
$$
\nu=\mu^-+\nu_0\ge 0,
$$
and obtain
\begin{equation}
\nu(\T) \le C \varkappa_s(F),             \label{v-estimate}
\end{equation}
and
\begin{multline*}
(\mathcal{P}*\nu)(z) \ge C_1 T(A(1-|z|)) \\
-C_2\Bigl(\log\frac 1 {1-|z|}\Bigr)^{s-1/2}\log\log \frac 1{1-|z|},\qquad z\in E^*_\eta.
\end{multline*}

Inequality (\ref{regt}) allows us to ignore the second summand in the right hand side
(replacing the constants if necessary):
\begin{equation}
(\mathcal{P}*\nu)(z)\ge C T(A(1-|z|)),\qquad z\in E^*_\eta.
\label{desiredestimate}
\end{equation}
\medskip

{\bf (F)} Finally, it turns out that $\mathcal{P}*\nu$ is too large on a too massive subset of the unit disc.

We  use the relation
\begin{equation}
\int_0\frac 1t\Bigl(\log\frac 1t\Bigr)^\beta\frac{dt}{T(At)}<\infty,
\label{b-convergence}
\end{equation}
which is just equivalent to \eqref{convergence}.

It follows from \eqref{v-estimate} and \eqref{desiredestimate} that
$$
m(E^*_\eta\cap r\mathbb{T})\le C \frac{\varkappa_s(F)}{T(A(1-r))}, \qquad 0<r<1,
$$
and hence, for every $r_1\in(0,1)$,
$$
\HH_\beta(E^*_\eta)\le \HH_\beta(r_1\D) +C\varkappa_s(F)
\int_0^{1-r_1}\frac 1\tau \Bigl(\log\frac 1\tau \Bigr)^{\beta}\frac{d\tau}{T(A\tau)}.
$$

Therefore,
$$
\limsup_{\varkappa_s(F)\to \infty}\frac{\HH_\beta(E^*_\eta)}{\varkappa_s(F)}\le
C \int_0^{1-r_1}\frac 1\tau \Bigl(\log\frac 1\tau \Bigr)^{\beta}
\frac{d\tau}{T(A\tau)},\qquad 0<r_1<1.
$$
By \eqref{b-convergence},
$$
\int_0^{1-r_1}\frac 1\tau \Bigl(\log \frac 1\tau \Bigr)^{\beta}
\frac{d\tau}{T(A\tau)}=o(1), \qquad r_1 \to 1,
$$
and we obtain
$$
\HH_\beta(E^*_\eta)=o(\varkappa_s(F)),\qquad \varkappa_s(F)\to \infty,
$$
that contradicts to \eqref{newbigarea}. This completes the proof of the theorem.
\medskip

\section{Proof of Theorem~\ref{thm2}}
\label{sect5}

{\bf (A)} Let us reformulate our problem.
Given $s\ge 1$, we are going to construct a finite singular measure $\mu$ on
$\mathbb{T}$, a set $\Omega\subset\D$ and a sequence of finite subsets $F_n$ of $\T$ such that
\begin{gather}
u(z)=(\mathcal{P}*\mu)(z)\ge M(1-|z|), \qquad z\in \Omega,\label{d12}\\
\inf_n \frac{\mathcal H_s(\Omega\cap {\mathcal K}_1(F_n))}{\varkappa_s(F_n)}>0,
\label{bigarea1}\\
\varkappa_s(F_n)\to\infty,\qquad n\to\infty,
\label{harmonicestimate}
\end{gather}
Then we verify that $\Omega\in\SU(\A_s)$;
if $f=\exp(-u-i\tilde{u})$,
where $\tilde{u}$ is the harmonic function conjugate to $u$, then
$$
\log |f(z)|\le -M(|z|),\qquad z\in\Omega.
$$

For $s<1$, we can just use $\Omega$ and $f$ constructed for $s=1$.
Therefore, from now on we assume that $s\ge 1$.

The measure $\mu$ and the set $\Omega$ will be constructed in a Cantor type
iterative procedure.
\medskip

{\bf (B)} Now we define several sequences of integer numbers which control
our iterative procedure. First, we set
$$
q_k=2^k,\qquad k\ge 1,
$$
and discretize our minorant $M$:
$$
m_1^0=1, \qquad m^0_k=\max\{1,\lfloor \lt M(1-2^{-q_k}) \rfloor  \}, \qquad k>1.
$$
where $\lfloor x \rfloor$ is the maximal integer number that does not exceed
$x$. We define $\mathcal W=\{k>1: m^0_{k-1}=m^0_k<m^0_{k+1}\}$.
Next, we modify a bit our sequence $\{m^0_k\}$. For $k>1$ we set $m^1_k=m^0_k+1$ if $k\in \mathcal W$,
and $m^1_k=m^0_k$ otherwise. The sequence $\{m^1_k\}$ does not differ too much from $\{m^0_k\}$, and for any
group $m^1_j=\ldots=m^1_k<m^1_{k+1}$ we have $m^1_{k+1}=m^1_k+1$.

Furthermore, we set
$$
m_1=1, \qquad m_k=\min\{m^1_k,2m_{k-1},m_{k-1}+s+1\}, \qquad k>1.
$$
The sequence $\{m_k\}$ is a lower regularization of $\{m^1_k\}$ of at most linear growth:
$m_k\le m_{k-1}+s+1$.

Finally, we set
$$
\mathcal Z=\{k\ge 1: m_k=m^0_k\}.
$$
Since
$$
\int_0 \frac 1t \Bigl(\log\frac 1 t \Bigr)^{s-1} \frac{dt}{M(1-t)}=\infty,
$$
we have
$$
\sum_{k\ge 1}q_k^s2^{-m_k}=\infty.
$$
Furthermore, since $m_k=m_{k-1}+s+1$ for large $k\not\in \mathcal Z$, we have
\begin{equation}
\sum_{1\le k\le N}q_k^s2^{-m_k}\asymp \sum_{k\in \mathcal Z:1\le k\le N}q_k^s2^{-m_k},
\qquad N\to\infty.
\label{d8}
\end{equation}
Thus, the set $\mathcal Z$ is massive enough to capture the growth asymptotics of $\sum q_k^s2^{-m_k}$.

Two more sequences $\{e_k\}$ and $\{L(k)\}$ are also defined in an inductive way.
Set $e_1=e_2=0$, $L(1)=1$, and, for $k>1$,

if $m_k>m_{k-1}$, then
\begin{equation}
L(k)=k, \qquad e_{k+1}=q_k-m_k,
\label{d1}
\end{equation}
otherwise, if $m_k=m_{k-1}$, then
\begin{equation*}
L(k)=L(k-1), \qquad e_{k+1}=e_k,
\label{d2}
\end{equation*}
Note that $e_k$ are non-negative,
\begin{equation}
e_{L(k)+1}=e_{k+1}, \qquad m_{L(k)}=m_k,\qquad k\ge 1.
\label{d3}
\end{equation}

Finally, we define a sequence $\{p_k\}$ of non-negative numbers by the relation
\begin{equation}
p_k=e_{k+1}-e_k,\qquad k\ge 1.
\label{d4}
\end{equation}

We choose $k_0$ such that $q_{k_0}-p_{k_0}-e_{k_0}=m_{k_0}>10$, $k_0=L(k_0)$, and
denote $\mathcal X=\{k\ge k_0:k=L(k)\}$, $\mathcal Y=\{k\ge k_0:k>L(k)\}$.
Clearly, $\mathcal Y\subset\mathcal Z$, and $L(x)\in \mathcal X$, $k\ge k_0$,
\begin{equation}
q_k=p_k+e_k+m_k,\qquad k\in\mathcal X.
\label{d16}
\end{equation}

Given two integer numbers $m\le n$ we denote by $[m,n]$
the set of all $k\in\NN$ such that $m\le k\le n$. If
$$
\mathcal Y=\cup_{l\ge 1}[v_l,w_l]
$$
is the (unique) representation such that
$$
w_l+1<v_{l+1},\qquad l\ge 1,
$$
then $w_l+1\in\mathcal X\cap\mathcal Z\cap\mathcal W$,
$m_{w_l+1}=m_{w_l}+1$, $l\ge 1$, and $m_k=m_{w_l}$, $k\in[v_l,w_l]$, $l\ge 1$.
Therefore, by \eqref{d8}, for $N\in\mathcal X$,
\begin{gather}
\sum_{k\in[1,N]}q_k^s2^{-m_k}\asymp \sum_{k\in \mathcal Z\cap[1,N]}q_k^s2^{-m_k}\notag\\
=\sum_{k\in \mathcal X\cap \mathcal Z\cap[1,N]}q_k^s2^{-m_k}
+\sum_{l:w_l<N}\sum_{k\in[v_l,w_l]}q_k^s2^{-m_k}\notag\\
\asymp\sum_{k\in \mathcal X\cap \mathcal Z\cap[1,N]}q_k^s2^{-m_k}
+\sum_{l:w_l<N}q_{w_l+1}^s2^{-m_{w_l}}\notag\\
\asymp\sum_{k\in \mathcal X\cap \mathcal Z\cap[1,N]}q_k^s2^{-m_k}
+\sum_{l:w_l<N}q_{w_l+1}^s2^{-m_{w_l+1}}\notag\\   \asymp
\sum_{k\in \mathcal X\cap \mathcal Z\cap[1,N]}q_k^s2^{-m_k},\qquad N\to\infty.
\label{d9}
\end{gather}
In particular, we conclude that the set $\mathcal X\cap\mathcal Z$ is infinite.
\medskip

{\bf (C)} Next, in a Cantor-type procedure, we define embedded families of intervals.
Their intersection will be the support of our measure $\mu$ to be produced later on.

For every $k\ge k_0$ we construct inductively a family of arcs on the unit circle
$\{J^{k}_{n,m}\}$, $n=1,2,\ldots,2^{e_k}$, $m=1,2,\ldots,2^{p_k}$, such that
$|J^{k}_{n,m}|= 2^{-q_{L(k)}}$.
For $k\in \mathcal X$, $k\ge k_0$, we will also construct families of subsets of the unit disc
$\Omega^{k-1}_{n,m}$, $n=1,2,\ldots,2^{e_k}$, $m=1,2,\ldots,2^{p_k}$, and finite subsets
$\mathfrak F^k$ of the unite circle.

We take the starting $2^{e_{k_0}+p_{k_0}}$ disjoint arcs of length $2^{-q_{k_0}}$
in an arbitrary way.
Now assuming that all $J^{k-1}_{*,*}$ are already constructed,
we define the arcs  $J^{k}_{n,m}$, $k$, $n=1,2,\ldots,2^{e_k}$, $m=1,2,\ldots,2^{p_k}$.

First, suppose that $k\in \mathcal X$. Then on each arc $J^{k-1}_{n,m}$ we take
a subarc $R^{k-1}_{n,m}$ which is the union of $2\cdot 2^{p_k}$ adjacent arcs
$J^{k}_{u,1},I^{k}_{u,1},J^{k}_{u,2},I^{k}_{u,2}$, \ldots, $J^{k}_{u,2^{p_k}},I^{k}_{u,2^{p_k}}$
of equal length $2^{-q_k}$, here $u=2^{p_{k-1}}n+m$. This is  possible because
$$
|J^{k-1}_{n,m}|=2^{-q_{L(k-1)}}\ge |R^{k-1}_{n,m}|.
$$
Indeed, $|R^{k-1}_{n,m}|= 2^{p_k-q_k+1}$, and by \eqref{d1}, \eqref{d3}, and \eqref{d16},
and the fact that $k,L(k-1)\in\mathcal X$, we have
\begin{gather*}
p_k-q_k+1+q_{L(k-1)}=-m_k-(e_{L(k-1)+1}-q_{L(k-1)})+1\\
=-m_k+m_{L(k-1)}+1\le 0.
\end{gather*}

Furthermore, we set
$$
\Omega^{k-1}_{n,m}=\bigl\{z:\arg z\in R^{k-1}_{n,m},\, 2^{-q_k}<1-|z|<2^{-q_{k-1}}\bigr\},
$$
see Figure 1. Finally, we define $\mathfrak F^k$
as the set of the middle points of the intervals $I^{k}_{*,*}$.

\begin{picture}(120,80)
\thicklines
\put(5,20){\line(1,0){50}}
\put(57,20){$\ldots$}
\put(65,20){\line(1,0){40}}

\put(10,32){\line(1,0){80}}
\put(10,72){\line(1,0){80}}
\put(10,32){\line(0,1){40}}
\put(90,32){\line(0,1){40}}

\thinlines
\put(10,20){\circle*{1}}
\put(100,20){\circle*{1}}
\put(90,20){\circle*{1}}
\put(10,15){\vector(1,0){80}}
\put(10,15){\vector(-1,0){0}}
\put(60,10){$R^{k-1}_{n,m}$}
\put(10,7){\vector(1,0){90}}
\put(10,7){\vector(-1,0){0}}
\put(80,2){$J^{k-1}_{n,m}$}
\put(45,-3){Figure 1.}
\put(10,20){\line(0,-1){15}}
\put(90,20){\line(0,-1){7}}
\put(100,20){\line(0,-1){15}}

\put(20,20){\circle*{1}}
\put(30,20){\circle*{1}}
\put(40,20){\circle*{1}}
\put(50,20){\circle*{1}}
\put(70,20){\circle*{1}}
\put(80,20){\circle*{1}}
\put(90,20){\circle*{1}}

\put(10,20){\line(0,1){6}}
\put(20,20){\line(0,1){6}}
\put(30,20){\line(0,1){6}}
\put(40,20){\line(0,1){6}}
\put(50,20){\line(0,1){6}}
\put(70,20){\line(0,1){6}}
\put(80,20){\line(0,1){6}}
\put(90,20){\line(0,1){6}}
\put(10,24){\vector(1,0){10}}
\put(10,24){\vector(-1,0){0}}
\put(20,24){\vector(1,0){10}}
\put(20,24){\vector(-1,0){0}}
\put(30,24){\vector(1,0){10}}
\put(30,24){\vector(-1,0){0}}
\put(40,24){\vector(1,0){10}}
\put(40,24){\vector(-1,0){0}}
\put(70,24){\vector(1,0){10}}
\put(70,24){\vector(-1,0){0}}
\put(80,24){\vector(1,0){10}}
\put(80,24){\vector(-1,0){0}}
\put(11,27){$J^{k}_{u,1}$}
\put(31,27){$J^{k}_{u,2}$}
\put(71,27){$J^{k}_{u,2^{p_k}}$}
\put(21,27){$I^{k}_{u,1}$}
\put(41,27){$I^{k}_{u,2}$}
\put(81,27){$I^{k}_{u,2^{p_k}}$}

\put(70,60){$\Omega^{k-1}_{n,m}$}

\put(90,32){\line(1,0){5}}
\put(90,72){\line(1,0){15}}

\put(93,20){\vector(0,1){12}}
\put(93,20){\vector(0,-1){0}}
\put(94,24){$2^{-q_k}$}
\put(102,20){\vector(0,1){52}}
\put(102,20){\vector(0,-1){0}}
\put(103,45){$2^{-q_{k-1}}$}

\end{picture}

\bigskip
\bigskip

In the case $k\in\mathcal Y$, we have $p_k=0$, $e_k+p_k=e_{k-1}+p_{k-1}$.
We define the set $\{J^{k}_{*,*}\}$ to be the same as
$\{J^{k-1}_{*,*}\}$ (the way the set is indexed depends on $e_k,p_k$ and
may be different from that for $J^{k-1}_{*,*}$).
\medskip

{\bf (D)} Next, for $k>k_0$ such that $k\in\mathcal X$, we verify the following statements:

\begin{enumerate}

\item
The sets $\Omega^{k-1}_{n,m}$ are disjoint;

\item
$\card\{J^k_{*,*}\}=2^{e_{k+1}}$;

\item
\begin{equation*}
\sum_{m,n}|J^k_{m,n}|=2^{-m_k};
\end{equation*}

\item We denote by $\mathcal{I}^k$ the set of all intervals of length
$2^{-q_k+1}$ located between two adjacent points of $\mathfrak F^k$.
Then
\begin{equation}
\sum_{I\in \mathcal{I}^k}|I|\Bigl(\log\frac1 {|I|}\Bigr)^s\asymp
q_k^s 2^{-m_k};
\label{partialentropy}
\end{equation}

\item We denote $\Omega^{k-1}=\cup_{m,n}\Omega^{k-1}_{m,n}$. Then
\begin{equation}
\mathcal H_s(\Omega^{k-1})\asymp q_k^s 2^{-m_k}.
\label{haka}
\end{equation}
\end{enumerate}

Statements (1) and (2) are straightforward.

Statement (3) follows from \eqref{d16}:
$$
\sum_{m,n}|J^k_{m,n}|=2^{-q_k+p_k+e_k}=2^{-m_k}.
$$

Statement (4) is an immediate consequence of statement (3).

A direct calculation gives
$$
\mathcal H_s(\Omega^{k-1}_{m,n})\asymp q_k^s|R^{k-1}_{m,n}|=q_k^s 2^{p_k-q_k+1}.
$$
By \eqref{d16},
$$
\mathcal H_s(\Omega^{k-1})\asymp q_k^s2^{p_k-q_k+e_k}=q_k^s2^{-m_k},
$$
and (5) follows.
\medskip

{\bf (E)} Let $\Omega=\cup_{k>k_0:k\in\mathcal X\cap \mathcal Z}\Omega^{k-1}$.
Here we verify \eqref{bigarea1} and \eqref{harmonicestimate}.
After that, to find an $s$-separated subset of $\Omega$ that belongs to $\U(\A^r_s)$
we just take a maximal $s$-separated (with a suitable constant $\e=\e(r)>0$)
subset $E$ of $\Omega$. Then, by Theorem~\ref{korteo}, the properties
\eqref{bigarea1} and \eqref{harmonicestimate} show that $E\in \U(\A^r_s)$.

Given $N>k_0$, denote
$$
\Omega_N=\cup_{k\in[k_0,N]\cap\mathcal X\cap \mathcal Z}\Omega^{k-1}.
$$
It follows from \eqref{d9} and \eqref{haka} that
$$
\mathcal H_s(\Omega_N)\asymp \sum_{k\le N}q_k^s2^{-m_k},\qquad N\to\infty.
$$

It remains to find $N_j\to\infty$, $j\to\infty$, and finite sets $F_{N_j}\subset\T$ such that
\begin{equation}
\mathcal{K}_1(F_{N_j})\supset\Omega_{N_j},
\label{inclusion}
\end{equation}
and
\begin{equation}
\varkappa_s(F_{N_j})\asymp\sum_{k\le N_j}q_k^s2^{-m_k},\qquad N_j\to\infty.
\label{entrestimate}
\end{equation}

We fix $N\in\mathcal X\cap\mathcal Z$, and set
\begin{equation*}
F_N=\bigcup_{k\in[k_0,N]\cap \mathcal X\cap\mathcal Z} \mathfrak F^k,
\end{equation*}
where $\mathfrak F^k$ are defined on step (C).

Relation (\ref{inclusion}) is now straightforward so it suffices to verify (\ref{entrestimate}).

To estimate the $s$-entropy of $F_N$, we use that $\mathcal{K}_1(F_{N})\supset\Omega_{N}$,
and that the groups of points in $\mathfrak F^k$, $k\in\mathcal X$, with distance
$2^{-q_{k}+1}$ between adjacent points, are situated
inside intervals of length $2^{-q_{L(k-1)}+1}$, with endpoints in $\mathfrak F^{L(k-1)}$.
It follows from \eqref{d9} and \eqref{partialentropy} that
$$
c\mathcal H_s(\Omega_N)\le
\varkappa_s(F_N)\le C\sum_{k\in [k_0,N]\cap \mathcal X\cap \mathcal Z} q_k^s 2^{-m_k}
\asymp \sum_{k\le N} q_k^s 2^{-m_k},\quad N\to\infty.
$$
This completes the proof of (\ref{entrestimate}).
\medskip

{\bf (F)} Finally, we define a measure $\nu$ on the unit circle (in a unique way)
by the following properties: (i) $\supp \nu\subset \cap_{k\ge k_0} \cup_{n,m} J^k_{n,m}$,
(ii) $\nu(J^{k}_{n,m})=2^{-e_{k+1}}$, $k\ge k_0$. Then for some sufficiently large $A$,
the measure $\mu=A\nu$ satisfies \eqref{d12}.
To verify this, we just use elementary estimates on the Poisson integral.
\medskip

\section{Approximation with restrictions in $C_A^\infty$.}
\label{sect6}

By duality, using a method of Havinson \cite{borichev5}, we can deduce
from Theorem~\ref{te2} a result on approximation by simple fractions
with restrictions on the coefficients in the space $C^\infty_A=C^\infty(\T)\cap H^\infty$.

The topology in $C_A^\infty$ is generated by the system of norms
$$
\|\phi\|_k = \|\phi^{(k)}\|_{L^2(\mathbb{T})},  \qquad  k\ge 0.
$$
If $C^n_A=C^n(\T)\cap H^\infty$, $n\ge 0$,
then convergence in $C_A^\infty$ means convergence in every $C^n_A$.

The dual to $C_A^\infty$ is $\A_1$ with the natural pairing: if $f\in \A_1$,
$\phi\in C_A^\infty$,
$$
f(z)=\sum_{k\ge 0} a_k z^k, \qquad
\phi(z)=\sum_0^\infty b_k z^k,
$$
then
$$
\langle\phi, f \rangle=\sum_{k\ge 0} a_k b_k.
$$
It is easily seen that, with such pairing, $(C^n_A)^*\subset \A^{n+2}_1$.

Furthermore, for $\phi_\lambda(z)=(1-z\lambda)^{-1}$, $|\lambda|<1$, we have
$$
\langle \phi_\lambda,f \rangle = f(\lambda).
$$
This immediately yields the following approximation result:
if $E\subset\D$ is a uniqueness set for $\A_1$, then the system
$\{(1-z\lambda)^{-1}\}_{\lambda \in E}$ is complete in $C_A^\infty$.

Here we prove that if $E\in\SU(\A_1)$, then it is possible to control
the coefficients of the approximating linear combinations.

\begin{thm} Let $E\in\SU(\A_1)$, and let a non-decreasing function
$M:(0,1)\to(0,\infty)$ satisfy the property
$$
\int_0 \frac 1t\frac{dt}{M(1-t)}<\infty.
$$
Then, for every $\phi\in C^\infty_A$, there exists a sequence of finite sums
$$
\omega_k(z)= \sum_{\lambda \in E} c_\lambda^{k}(1-z\lambda)^{-1},
$$
such that
$$
\omega_k\to\phi \mbox{\ in $C^\infty_A$},\qquad k\to\infty,
$$
and
$$
\sum_{\lambda \in E}|c^k_\lambda|e^{-M(|\lambda|)}\to 0, \qquad k\to\infty.
$$
\label{tha}
\end{thm}

To obtain this result we use the following lemma.

\begin{lem} Given $n\in\NN$, a sequence $\{\lambda_m\}_{m\ge 1}$ of points in $\D$
and a sequence of positive numbers $\{t_m\}_{m\ge 1}$, the following statements are equivalent:
\begin{enumerate}
\item For every $\phi\in C^n_A$ there exist finite sums
$$
\omega_k(z)= \sum_{m\ge 1} c_m^{k}(1-z\lambda_m)^{-1},
$$
such that
$$
\omega_k\to\phi \mbox{\ in $C^n_A$},\qquad k\to\infty,
$$
and such that
$$
\sum_{m\ge 1} \frac{|c^k_m|}{t_m}\to 0, \qquad k\to\infty.
$$
\item If $f\in (C^n_A)^*$ and if
$$
\Bigl|\sum_{m\ge 1}c_m f(\lambda_m)\Bigr| \le  \sup_{m\ge 1} \frac{|c_m|}{t_m}
$$
for any finite sequence $\{c_m\}$, then $f=0$.
\end{enumerate}
\label{lele}
\end{lem}

We can just choose in this Lemma $t_m=\exp{M(|\lambda_m|)}$, and apply
Theorem~\ref{thm1} and relation \eqref{regt}
to obtain Theorem~\ref{tha}.

Lemma~\ref{lele} follows from Theorem~1 in \cite{borichev1} if we set
$X=C^n_A$, $p(\phi)=\|\phi\|_{C^n_A}$, and $p_1(\{c_m\})= \sup_{m\ge 1} \frac{|c_m|}{t_m}$
(as in Theorem~2 of \cite{borichev1}).
\medskip

\section{Final remarks.}
\label{sect7}

We could compare our results to those by Pau and Thomas \cite{borichev4} concerning
the minorants for pairs $(H^\infty,\E)$, for some special classes $\E\subset\SU(H^\infty)$.

Given a subset $E$ of the unit disc, they define a function
$\phi_E:\T\to\Z_+\cup\{+\infty\}$ by
$$
\phi_E(\zeta)=\card\bigl(E\cap\Gamma(\zeta)\bigr),
$$
where $\Gamma(\zeta)$ is the Stolz angle (of fixed aperture) with vertex at $\zeta\in\T$.
Furthermore, $h_E:\Z_+\to[0,1]$ is defined by
$$
h_E(n)=m\bigl(\zeta\in\T:\phi_E(\zeta)\ge n\bigr).
$$
Given a non-increasing sequence of positive numbers $w=\{w_n\}_{n\ge 1}$, they consider
the class $\E_w$ of hyperbolically separated subsets $E$ of the unit disc such that
$$
\sum_{n\ge 1}h_E(n)w_n=\infty.
$$
Clearly, $\E_w\ne\emptyset$ if and only if $\sum_{n\ge 1}w_n=\infty$.

\begin{thm} {\rm(Pau--Thomas, \cite{borichev4})} Let $w$ be fixed,
and let $M:(0,1)\to(0,\infty)$ be a non-decreasing function.
For every $f\in H^\infty$, $E\in\E_w$,
the estimate
\begin{equation}
\log|f(z)|\le -M(|z|),\qquad z\in E
\label{d10}
\end{equation}
implies that $f=0$ if and only if
\begin{equation}
\sum_{n\ge 1}\frac{w_n}{M(1-2^{-n})}<\infty.
\label{d11}
\end{equation}
\end{thm}

Here we show how to deduce the necessity part of this result using
our constructions. Pau and Thomas use a different method in \cite{borichev4}.

Suppose that \eqref{d11} does not hold. Then
$$
\sum_{n\ge 1}\frac{1}{M(1-2^{-n})}=\infty,
$$
and hence, condition \eqref{divergence} holds. Using the construction
in the proof of Theorem~\ref{thm2} (with $s=1$),
we can find $f\in H^\infty\setminus\{0\}$ and a hyperbolically separated $E\subset\D$
such that \eqref{d10} is satisfied, and
$$
h_E(2q_k)\ge c\cdot 2^{-m_k}, \qquad k\in\mathcal X\cap\mathcal Z.
$$
Since $h_E$ is monotonic, this estimate extends from $k=w_l+1$ to all $k\in[v_l,w_l],$
and hence, we have the same estimate for all $k\in\mathcal Z$.
Therefore, by monotonicity of $\{w_n\}$, we obtain that
\begin{gather*}
+\infty=\sum_{n\ge 1}\frac{w_n}{M(1-2^{-n})}\le
c+c\sum_{k\ge 1}\Bigl(\sum_{q_k\le n<q_{k+1}}w_n\Bigr)2^{-m_k}\\ \le
c+c\sum_{k\in \mathcal Z}\Bigl(\sum_{q_k\le n<q_{k+1}}w_n\Bigr)2^{-m_k}\le
c+c\sum_{k\in \mathcal Z}\sum_{q_k\le n<q_{k+1}} w_n h_E(n)\\ \le
c+c\sum_{k\ge 1} w_n h_E(n),
\end{gather*}
and we conclude that $E\in\E_w$.

Next we present an example showing how different are
the classes $\SU(\A_s)$ for $s\le 1$ and for $s>1$.

\begin{exemp} Given $s>1$, we produce here a set in
$\SU(\A_s)$, which contains no hyperbolically separated non-Blaschke sequence.

We choose two increasing sequences of integer numbers $\{n_k\}_{k\ge 0}$, $\{m_k\}_{k\ge 0}$,
such that
$$
m_kn_k^s 2^{-n_k/2} \to \infty,\qquad k\to\infty,
$$
and
$$
\sum_{k\ge 0}m_kn_k 2^{-n_k/2} < 1.
$$
After that, we find a family $\{\phi_{k,j}\}_{k\ge 0, 1\le j\le m_k}$ such that
the intervals
$$
I_{k,j}= \{e^{i\theta}:\phi_{k,j}\le \theta\le \phi_{k,j}+2^{-n_k/2}\}\subset (0,1)
$$
are pairwise disjoint and consider the sets $E_k,E\subset\D$ and $F_k\subset\T$, $k\ge 0$,
defined as follows
\begin{gather*}
E_{k,j}=\{re^{i\theta}: \theta \in I_{k,j},\, 2^{-n_k}<1-r<2^{-n_k/2}\},\\
E_k=\cup_{j=1}^{m_k} E_{k,j}, \qquad E=\cup_{k\ge 0} E_k,\\
F_{k,j}=\{\exp(2\pi i l 2^{-n_k})\}_{l=1}^{2^{n_k}}\cap I_{k,j}, \qquad
F_k=\cup_{j=1}^{m_k} F_{k,j}.
\end{gather*}

It follows from the construction that all $E_{k,j}$ are disjoint,
$E_k\subset \mathcal{K}_1(F_k)$, and
$$
\mathcal H_s(E_k)\asymp \varkappa_s(F_k) \asymp m_kn_k^s 2^{-n_k/2} \to \infty,
\qquad k\to\infty.
$$
Therefore, by Theorem~\ref{korteo}, we conclude that $E\in\SU(\A_s)$.

On the other hand,
$$
\mathcal H_1(E)\asymp \sum_{k\ge 0} m_kn_k 2^{-n_k/2} < \infty,
$$
which yields that $E$ contains no hyperbolically separated
non-Blaschke sequence.
\end{exemp}

\bigskip
\bigskip

\noindent \textsc{Alexander Borichev, Department of Mathematics,}

\noindent \textsc{University of Bordeaux I, 351, cours de la Lib\'eration, 33405 Talence, France}

\noindent\textsl{E-mail}: \texttt{borichev@math.u-bordeaux.fr}
\medskip

\noindent \textsc{Yurii Lyubarskii, Department of Mathematical
Sciences, Norwegian University of Science and Technology,
N-7491 Trondheim, Norway}

\noindent\textsl{E-mail}: \texttt{yura@math.ntnu.no}

\end{document}